\documentclass[11pt]{article}%
\usepackage{graphicx}
\usepackage{subcaption}
\usepackage{makeidx}
\usepackage{amssymb}
\usepackage[all]{xy}
\usepackage{float}
\usepackage[latin1]{inputenc}
\usepackage{amsfonts}
\usepackage{amsmath}
\usepackage{amssymb}
\usepackage{url}
\usepackage{hyperref}
\usepackage{graphicx}%
\usepackage{mathabx}
\usepackage{wrapfig}

\usepackage{amsthm}
\usepackage{url}
\usepackage{hyperref}
\providecommand{\U}[1]{\protect\rule{.1in}{.1in}}

\usepackage{xcolor}

\newtheorem{prop}{Proposition}[section]

\newtheorem{theo}[prop]{Theorem}

\newcommand{\vertiii}[1]{{\left\vert\kern-0.25ex\left\vert\kern-0.25ex\left\vert #1 
    \right\vert\kern-0.25ex\right\vert\kern-0.25ex\right\vert}}

\def\tr{\mbox{\rm Tr}}

\newcommand{\LL}{\mathbb{L}}

\newcommand{\La}{ {\cal L }}

\newcommand{\Sa}{ {\cal S }}

\newcommand{\Oa}{ {\cal O }}

\newcommand{\cqfd}{\hfill\blbx \\}
\def\blbx{\hbox{\vrule height 5pt width 5pt depth 0pt}\medskip}

\def \LL{\mathbb{L}}

\makeatletter
\DeclareRobustCommand\frownotimes{\mathbin{\mathpalette\frown@otimes\relax}}
\newcommand{\frown@otimes}[2]{%
  \vbox{
    \ialign{##\cr
      \hidewidth$\m@th#1{}_\frown$\kern-\scriptspace\hidewidth\cr
      \noalign{\nointerlineskip\kern-.1pt}
      $\m@th#1\otimes$\cr
    }%
  }%
}
\makeatother

\begin{document}
\date{}

  \title{A Taylor expansion of the square root matrix  function}
  \author{P. Del Moral\footnote{INRIA Bordeaux Sud Ouest Center, 33 405 Talence, pierre.del-moral@inria.fr}, ~A. Niclas\footnote{ENS Lyon, 46 All\'ee d'Italie, 69007 Lyon, angele.niclas@ens-lyon.fr}}


\maketitle

\begin{abstract}
This short note provides an explicit description of the Fr\'echet derivatives of the principal 
square root matrix  function at any order. We present an original formulation that
allows to compute sequentially the Fr\'echet derivatives of the matrix square root at any order
starting from the first order derivative. 
 A Taylor expansion at any order with an integral remainder term is also provided,
yielding the first result of this type for this class of matrix  function.



\end{abstract}


\section{Introduction}

The computation of matrix square roots arise in a variety of application domains, 
including in physics, signal processing,
optimal control theory, and many others. 
The literature abounds with numerical 
techniques for computing matrix square roots, see for instance~\cite{almohy,davies,deadman,higham2,higham3,higham4,meini}.
Perturbative techniques often resume to Lipschitz type estimates~\cite{hemmen} or on the refined analysis of the first order Fr\'echet derivative of the principal 
square root matrix function; see for instance~\cite{almohy,almohy2,almohy3,schmitt}, as well as chapter X in the seminal book by R. Bhatia~\cite{Bhatia} and references therein. We also refer to the article~\cite{Cardoso}
for a first order analysis of more general matrix $n$-th roots. For further details on the $n$-th roots of matrices we refer to~\cite{psarrakos}.

The purpose of this article is to derive an explicit description of the Fr\'echet derivatives of the principal 
square root matrix  function at any order. We also provide a non asymptotic 
Taylor expansion at any order, with computable estimates of the integral remainder terms.
 These expansions provide a perturbation computation of the square root $\sqrt{A+H}$ of a positive definite matrix $A$ perturbed by some symmetric matrix $H$, as soon as $A+\epsilon~H$ is positive semidefinite for any $\epsilon\in [0,1]$. 
 
 We underline that the perturbation analysis developed in this article differs from Taylor expansion type techniques often used to
 define functions on the spectrum of diagonalizable matrices via Jordan canonical forms. This Sylvester's formulation of matrices are related to the Sylvester matrix theorem (a.k.a. Lagrange-Sylvester interpolation)
 which allows to express an analytic function of a matrix in terms of its eigenvalues and eigenvectors.
 For a more thorough discussion on these interpolation techniques we refer to the first
 chapter  in the seminal book by N. J. Higham~\cite{higham}.

This study has been motivated by applications in signal processing and more particularly in the 
analysis of Ensemble Kalman-Bucy filters~\cite{dbn17}. In this context, the sample 
and interacting covariance matrices satisfy a stochastic matrix Riccati diffusion. The diffusion term
depends on the matrix square root of the sample covariance. The perturbation analysis developed in this article is used
to derive non asymptotic Taylor-type expansions of stochastic matrix Riccati flows  w.r.t. some perturbation parameter.

We denote by $\Sa_{r}$ the space of symmetric $(r\times r)$-matrices $A$ with real entries
equipped with  the $\LL_2$-norm $\Vert A \Vert=\Vert A \Vert_2=\sqrt{\lambda_{ max}(A^2)}$ or the Frobenius norm  $\Vert A \Vert=\Vert A \Vert_F=\sqrt{\tr(A^2)}$. We recall that these norms are 
equivalent and $\Vert A \Vert_2\leq \Vert A \Vert_F\leq \sqrt{r}~\Vert A \Vert_2$.

We let
 $ \Sa_r^0\subset\Sa_r$  { be the closed convex cone} of  {positive} semi-definite matrices,  and 
its interior $\Sa_r^+\subset \Sa_r^0$  which resumes to {the open subset} of  {positive definite}  matrices. 
We consider the principal square root  function 
$$
\varphi~:~Q\in \Sa^+_r\mapsto
\varphi(Q)=Q^{1/2}\in \Sa_r^+
$$ 
For any $Q_1,Q_2\in\Sa_{r}^+$ we have the Ando-Hemmen inequality
\begin{equation}\label{square-root-key-estimate}
\Vert \varphi(Q_1)- \varphi(Q_2)\Vert \leq \left[\lambda^{1/2}_{min}(Q_1)+\lambda^{1/2}_{min}(Q_2)\right]^{-1}~\Vert Q_1- Q_2\Vert
\end{equation}
for any unitary invariant matrix norm $\Vert . \Vert$. See for instance Theorem 6.2 on page 135 in~\cite{higham}, as well as Proposition 3.2  in~\cite{hemmen}. For a more thorough discussion on the geometric properties of positive semidefinite matrices and square roots we refer to~\cite{hiriart}.

We let $\La(\Sa_r,\Sa_{r})$ be the set of bounded linear  functions from $\Sa_r$ into itself. Let  {$\Oa_r\subset\Sa_r$ be a non empty open and convex subset of $\Sa_r$}.
We recall that a mapping  {$\Upsilon:\Oa_{r}\mapsto \Sa_r$} defined in some domain $\Oa_{r}$
is Fr\'echet differentiable at some $A\in \Oa_r$ if there exists a 
continuous linear  function
$
\nabla \Upsilon(A)\in  \La(\Sa_r,\Sa_{r})
$
such that
$$
\lim_{\Vert H\Vert\rightarrow 0}\Vert H\Vert^{-1}
\Vert \Upsilon(A+H)-\Upsilon(A)-\nabla \Upsilon(A)\cdot H \Vert=0
$$
In other words, for any given $A\in \Oa_r$ and $\epsilon>0$ there exists some $\delta>0$ such that 
$$
\Vert H\Vert\leq \delta\Longrightarrow A+H\in \Oa_r\quad\mbox{\rm and}\quad
\Vert \Upsilon(A+H)-\Upsilon(A)-\nabla \Upsilon(A)\cdot H \Vert\leq \epsilon~\Vert H\Vert
$$
The l.h.s. condition is met for $\Oa_r=\Sa^+_r\ni A$. We check this claim using Weyl's inequality 
$$
\lambda_{ min}(A+H)\geq \lambda_{ min}(A)+\lambda_{ min}(H)\geq  \lambda_{ min}(A)-\Vert H\Vert_2
$$
This shows that
$$
\Vert H\Vert_2< \lambda_{ min}(A)\Longrightarrow A+H\in \Sa^+_r
$$
 {The function $\Upsilon$ is said to be Fr\'echet  differentiable on $\Oa_r$ when the mapping $$\nabla \Upsilon~:~A\in \Oa_r
\mapsto \nabla \Upsilon(A)\in  \La(\Sa_r,\Sa_{r})$$ is continuous. 
Higher Fr\'echet derivatives are defined in a similar way. For instance, the mapping $\Upsilon$  is twice Fr\'echet differentiable at $A\in \Oa_r$ when
the mapping
$
\nabla \Upsilon
$
is also Fr\'echet  differentiable at $A\in \Oa_r$}. Identifying $ \La(\Sa_r,\La(\Sa_r,\Sa_{r}))$ with the set 
$ \La(\Sa_r\times\Sa_r,\Sa_{r})$ of continuous bilinear maps from $(\Sa_r\times\Sa_r)$ into $\Sa_r$, the second derivative
$$\nabla^2 \Upsilon~:~A \in \Oa_r \mapsto \nabla^2 \Upsilon(A)\in  \La(\Sa_r\times\Sa_r,\Sa_{r})$$ is defined by a continuous and symmetric billinear map $\nabla^2 \Upsilon(A)$ such that the limit
$$
\lim_{\Vert H_2\Vert\rightarrow 0}\Vert H_2\Vert^{-1}
\Vert \nabla \Upsilon(A+H_2)\cdot H_1-\nabla \Upsilon(A)\cdot H_1-\nabla^2 \Upsilon(A)\cdot (H_1,H_2) \Vert=0
$$ 
exists uniformly w.r.t. $H_1\in\Sa_r$ in bounded sets. The polarization formula 
$$
\nabla^2 \Upsilon(A)\cdot (H_1,H_2)=\frac{1}{4}~\left[
\nabla^2 \Upsilon(A)\cdot (H_1+H_2,H_1+H_2)-\nabla^2 \Upsilon(A)\cdot (H_1-H_2,H_1-H_2)\right]
$$
shows that it suffices to compute the second order derivatives $\nabla^2 \Upsilon(A)\cdot (H,H)$ in the same direction $H=H_1=H_2$. Identifying
$\La(\Sa_r\times\Sa_r,\Sa_{r})$ with $\La(\Sa_r\otimes\Sa_r,\Sa_{r})$ sometimes we set $\nabla^2 \Upsilon(A)\cdot H^{\otimes 2}$ instead of $\nabla^2 \Upsilon(A)\cdot (H,H)$. For a more detailed discussion on these tensor product identifications of symmetric multilinear maps we refer to chapter 5 in~\cite{dudley}. 

Higher Fr\'echet derivatives $\nabla^n \Upsilon(A)\cdot H^{\otimes n}$ of order $n$ are defined recursively in a similar manner. For a more thorough discussion on higher Fr\'echet derivatives, we refer
the reader to the seminal books of Cartan~\cite{cartan} and Dieudonn\'e~\cite{dieudonne}, 
as well as section 5.2 in the book by Dudley and Norvaisa~\cite{dudley} and the
article by Higham and Relton~\cite{higham-relton-2014}. The latter addresses the general case, as well as the matrix exponential and inverse.
 
In the further development of this article, by symmetry arguments we only consider differentials in a given direction $H$. 
To simplify the
presentation, sometimes we write $\nabla^n\Upsilon(A)\cdot H$ instead of $\nabla^n\Upsilon(A)\cdot H^{\otimes n}$. The $n$-th 
derivatives 
$\nabla^n \Upsilon(A)\cdot (H_1,\ldots,H_n)$ in $n$ different directions are defined as above by
polarization of $n$-linear symmetric operators, see for instance theorem 5.6 in~\cite{dudley}. The symmetry property of the $n$-linear mappings $\nabla^n \Upsilon(A)$ is a consequence of Schwarz theorem (see for instance theorem 5.27 in~~\cite{dudley}).

 Let  $\Upsilon:\Oa_{r}\mapsto \Sa_r$ be a Fr\'echet differentiable mapping at any order at some $A\in \Oa_r$. Given some $H\in \Sa_r$ s.t. $A+H$ is included in $\Oa_{r}$ we have
 \begin{equation}\label{Taylor-remainder}
\Upsilon(A+H)=\Upsilon(A)+
\sum_{1\leq k\leq n}~\frac{1}{k!}~\nabla^k\Upsilon(A)\cdot H+\overline{\nabla}^{n+1}\Upsilon\left[A,H\right]
\end{equation}
with the $(n+1)$-th order remainder
 function in the Taylor expansion given
\begin{eqnarray*}
\overline{\nabla}^{n+1}\Upsilon\left[A,H\right]&:=&\frac{1}{n!}~\int_0^1~(1-\epsilon)^{n}~\nabla^{n+1}\Upsilon\left(A+\epsilon~H\right)\cdot H~d\epsilon
\end{eqnarray*}
Using the convexity of the set $\Oa_r$, we underline that the line segment $\epsilon\in [0,1]\mapsto A+\epsilon H$ joining $A\in\Oa_r$ to $A+H\in \Oa_r$
is included in $\Oa_{r}$; that is, we have that
\begin{equation}\label{ref-convex}
A\in \Oa_r\quad\mbox{\rm and}\quad
B=A+H\in  \Oa_+\Longrightarrow\forall\epsilon\in [0,1]\qquad A+\epsilon H=(1-\epsilon)~A+\epsilon~B\in \Oa_r 
\end{equation}

For a more detailed account on Taylor's formulae with integral remainders for smooth functions on open convex subsets of Banach spaces with values in another Banach space we refer the reader to section 5.3 in the book by Dudley and Norvaisa~\cite{dudley}.

We also consider the multi-linear operator norm
$$
\vertiii{\nabla^n\Upsilon(P)}=\sup_{\Vert H\Vert=1} 
\Vert \nabla^n\Upsilon(P)\cdot H\Vert$$
In the further development of this  article $C_{n}=\displaystyle\frac{1}{n+1}~\left(\begin{array}{c}
2n\\
n
\end{array}
\right)$ stands for the Catalan number.

\begin{theo}\label{lem-square-root-taylor}
The square root  function $\varphi ~:~Q\in \Sa^+_r\mapsto
\varphi(Q)=Q^{1/2}\in \Sa_r^+$ is Fr\'echet differentiable at any order on $\Sa_r^+$ with the first order derivative given
for any
$(A,H)\in (\Sa_r^+ \times \Sa_r)$ 
 by the formula
\begin{eqnarray}\label{induction-square-root-n-nabla}
\nabla \varphi(A)\cdot H&=&\int_0^{\infty}e^{-t\varphi(A)}~H~e^{-t\varphi(A)}~dt
\end{eqnarray}
The higher order derivatives are defined inductively for any $n\geq 2$ by the formula
 \begin{equation}\label{induction-square-root-n-nabla-bis}
\begin{array}{l}
\displaystyle \nabla^n \varphi(A)\cdot H\\
\\
=-\displaystyle\nabla \varphi(A)\cdot \left[
\sum_{p+q=n-2}~\frac{n!}{(p+1)!(q+1)!}~
\left[
\nabla^{p+1} \varphi(A)\cdot H\right]~\left[\nabla^{q+1} \varphi(A)\cdot H\right]\right]
\end{array}
\end{equation}
In the above display, the summation is taken over all integers $p,q\geq 0$ s.t. $p+q=n-2$. 
Assume that $A$ and $A+H\in \Sa_r^+$. In this situation the function $\varphi$ has a Taylor expansion (\ref{Taylor-remainder}) at any order.
In addition, for any $n\geq 0$ 
we have
the estimates
\begin{eqnarray}
\vertiii{
\nabla^{n+1} \varphi(A)}&\leq& K^{n}~(n+1)!~C_n~2^{-(2n+1)}~\lambda_{ min}(A)^{-(n+1/2)}\nonumber\\
\Vert \overline{\nabla}^{n+1}\varphi\, [A,H]\Vert&\leq &K^{n}~(n+1)~C_n~2^{-2n}~\lambda_{ min}(A)^{-(n+1/2)}
~\Vert H\Vert^{n+1}\label{estimate-square-root-n-nabla}
\end{eqnarray}
where $K=\sqrt{r}$ for the Frobenius norm, and $K=1$ for the $\LL_2$-norm. 

\end{theo}

We end this section with some comments on the above theorem. 

Firstly, arguing as in (\ref{ref-convex})  the convexity of the set $\Sa_r^+$ ensures that
the line segment joining the matrix $A\in \Sa_r^+$ to any matrix $B=A+H\in \Sa_r^+$ is always included in $\Sa_r^+$. The terminal state condition $B=A+H\in \Sa_r^+$ is  met for any  $H\in \Sa_r$ s.t. $\lambda_{ min}(A)>0\vee(-\lambda_{ min}(H))$. This condition is also clearly met for any $H\in\Sa_r^0$.

 The inductive formula (\ref{induction-square-root-n-nabla-bis}) allows to compute sequentially the Fr\'echet derivatives of the matrix square root at any order
starting from the first order derivative. For instance the second Fr\'echet derivative is given by
$$
\nabla^2 \varphi(A)\cdot H=-2\nabla \varphi(A)\cdot \left[
\nabla \varphi(A)\cdot H\right]^2
$$
In this situation, using (\ref{estimate-square-root-n-nabla}) for any $A\in\Sa_r^+$ and $B=A+H\in \Sa^+_r$ we find that
$$
\Vert \varphi(B)-\varphi(A)-\nabla \varphi(A)\cdot (B-A)+\nabla \varphi(A)\cdot \left[
\nabla \varphi(A)\cdot (B-A)\right]^2\Vert_2\leq \frac{3}{8}~\lambda_{ min}(A)^{-5/2}
~\Vert B-A\Vert^{3}_2
$$

As mentioned in the introduction
several alternative representations of the Fr\'echet derivative of the square root  function can be found in the literature. To better connect our work with existing results we end this section around this theme.

As shown in~\cite{bhatia2}, the integral representation
of the square root matrix  function is given in terms of the resolvent of $-A$ by the formula
$$
\begin{array}{l}
\displaystyle\varphi(A)=\frac{1}{\pi}~\int_0^{\infty}A(tI+A)^{-1}t^{-1/2}dt\\
\\
\displaystyle\Longrightarrow \nabla \varphi(A)\cdot H=\frac{1}{\pi}~\int_0^{\infty} (tI+A)^{-1}~H~(tI+A)^{-1}t^{1/2}dt
\end{array}
$$
The last assertion is proved using a simple differentiation under the integral sign (invoking the dominated
convergence theorem). The article~\cite{Cardoso} also extends this integral formulae to more general
  $n$-th roots matrix functions. The article~\cite{bhatia2} (see formula (15)) also provides an 
  alternative formulation in terms of the exponential matrix of $A$; namely
  $$
  \nabla \varphi(A)\cdot H=\frac{1}{2\sqrt{\pi}}~\int_0^{\infty}\left[\int_0^t~e^{-sA}~H~e^{-(t-s) A}~ds\right]~t^{-3/2}dt 
  $$

 It is well know that
the Fr\'echet derivative  $X=\nabla \varphi(A)\cdot H$ given in (\ref{induction-square-root-n-nabla})
is the unique solution of the Sylvester equation~\cite{sylvester} given by
$$
\psi(A)=A^2\Longrightarrow \nabla\psi(\varphi(A))\cdot X=\varphi(A)~X+X~\varphi(A)=H\quad\mbox{\rm and}\quad (\nabla \varphi)(\psi(A))=\left[\nabla\psi(A)\right]^{-1}
$$
See for instance, section 6.1 in~\cite{higham}, and the article~\cite{bhatia2}. 
The Sylvester equation stated above is a particular case of the algebraic Riccati equation. It can also be regarded as
a Lyapunov equation. In this connection there are no surprise that  (\ref{induction-square-root-n-nabla}) coincides
with the rather well know solution of the continuous Lyapunov equation. This integral formulation is closely related to the notion of controllability Gramian of a linear dynamical system with drift matrix $-\varphi(A)$, see for instance~\cite{brockett}.

The literature also abounds with numerical techniques for solving of the
Sylvester equation, see for instance the recent review by V. Simoncini~\cite{simoncini} and references therein.

The formulae (\ref{induction-square-root-n-nabla-bis}) for higher terms in the Taylor series for the square root
provide a polynomial-type perturbation approximation of the square root at any order. These non asymptotic expansions have been used in~\cite{dbn17,dbn17-2} to analyze the fluctuation as well as  the bias of the square root function of Wishart matrices and sample covariance matrices associated with stochastic Riccati equations arising in Ensemble-Kalman-Bucy filter theory.

\section{Proof of theorem~\ref{lem-square-root-taylor}}
Any (symmetric) 
square roots $\varphi(A)$ and $\varphi(B)$ of matrices $A,B\in\Sa_r^+$ satisfy the Sylvester equation
$$
\varphi(A)\left(\varphi(A)-\varphi(B)\right)+\left(\varphi(A)-\varphi(B)\right)\varphi(A)=(A-B)+\left(\varphi(A)-\varphi(B)\right)^2:=C
$$
When $\varphi(A)>0$ we have
$$
\begin{array}{l}
\displaystyle
Z_t:=e^{-t\varphi(A)}Ce^{-t\varphi(A)}~\longrightarrow_{t\rightarrow\infty}~ 0\\
\\
\displaystyle\Longrightarrow~-
\partial_tZ_t=\varphi(A)~Z_t+Z_t~\varphi(A)\\
\\
\displaystyle\Longrightarrow Z_0=C=\varphi(A)~\left[\int_0^{\infty}Z_t~dt\right]+\left[\int_0^{\infty}Z_t~dt\right]~\varphi(A)
\end{array}$$
This implies that
\begin{equation}\label{square-root-first-order}
\left(\varphi(B)-\varphi(A)\right)=\int_0^{\infty}e^{-t\varphi(A)}~(B-A)~e^{-t\varphi(A)}~dt
-\int_0^{\infty}e^{-t\varphi(A)}\left(\varphi(B)-\varphi(A)\right)^2e^{-t\varphi(A)}~dt
\end{equation}
We set $B:=A+H$. Under our assumptions $B\in\Sa_r^+$. 
Using (\ref{square-root-key-estimate}) we conclude that
\begin{equation}\label{induction-square-root-integral-pre}
\nabla \varphi(A)\cdot H=\int_0^{\infty}e^{-t\varphi(A)}~H~e^{-t\varphi(A)}~dt\Longrightarrow\vertiii{\nabla \varphi(A)}\leq ~K~2^{-1}~\lambda_{ min}^{-1/2}(A)
\end{equation}
We check the last assertion using the fact that
$$
\varphi(A)>0\Longrightarrow-\lambda_{ min}(\varphi(A))=\rho(-\varphi(A))=-\lambda_{ min}^{1/2}(A)
$$
This yields the remainder formula
\begin{equation}\label{induction-square-root-integral}
\begin{array}{l}
\overline{\nabla}^{2}\varphi\left[A,B-A\right]=\varphi(B)-\varphi(A)-\nabla \varphi(A)\cdot (B-A)=-\nabla \varphi(A)\cdot \left(\varphi(B)-\varphi(A)\right)^2\\
\\
\Longrightarrow \Vert\overline{\nabla}^{2}\varphi\left[A,B-A\right]\Vert\leq \vertiii{\nabla \varphi(A)}~\Vert\varphi(B)-\varphi(A)\Vert^2\leq K~2^{-1}~\lambda_{ min}^{-3/2}(A)~\Vert A-B\Vert^2
\end{array}
\end{equation}
The last assertion is a consequence of the Ando-Hemmen inequality (\ref{square-root-key-estimate}) and the estimate (\ref{induction-square-root-integral-pre}).
This ends the proof of the Taylor expansion at rank $n=1$.
We set  
$$
T_n(A,H)=\sum_{1\leq k\leq n}\frac{1}{k!}~
\displaystyle\partial^k \varphi(A)\cdot H
$$
with the collection of matrices $\partial^k \varphi(A)\cdot H$ defined by 
\begin{equation}\label{partial1-definition}
\partial^1 \varphi(A)\cdot H=\nabla \varphi(A)\cdot H
\end{equation}
 and for any $n\geq 2$ by
the induction
\begin{equation}\label{induction-square-root-n}
\begin{array}{l}
\displaystyle \partial^n \varphi(A)\cdot H\\
\\
:=-\displaystyle\nabla \varphi(A)\cdot \left[
\sum_{p+q=n-2}~\frac{n!}{(p+1)!(q+1)!}~
\left[
\partial^{p+1} \varphi(A)\cdot H\right]~\left[\partial^{q+1} \varphi(A)\cdot H\right]\right]
\end{array}
\end{equation}
In the above display, the summation is taken over all integers $p,q\geq 0$ s.t. $p+q=n-2$. 
We prove (\ref{estimate-square-root-n-nabla}) by induction on the parameter $n$. First, we prove that
\begin{equation}\label{catalan-definition}
\forall 1\leq k\leq n\qquad
 \frac{1}{k!}~\vertiii{
\partial^{k} \varphi(A)}\leq  \frac{K^{k-1}}{2\lambda_{ min}^{{1}/{2}}(A)}
~\frac{C_{k-1}}{2^{2(k-1)}\lambda_{ min}(A)^{k-1}}~
\end{equation}
By (\ref{induction-square-root-integral-pre}) and (\ref{partial1-definition}) this assertion is clearly met for $n=1$. Assume that the above estimates (\ref{catalan-definition}) are met for any $1\leq k<n$.
Combining (\ref{induction-square-root-integral-pre}) with (\ref{induction-square-root-n}) we find that
$$
\begin{array}{l}
\displaystyle\frac{1}{n!}~\displaystyle \vertiii{\partial^{n} \varphi(A)}
\leq  \frac{K}{2\lambda_{ min}^{{1}/{2}}(A)} \left[
\sum_{p+q=n,~p,q\geq 1}~\frac{1}{p!}~
\vertiii{
\partial^{p} \varphi(A)}~\frac{1}{q!}~\vertiii{\partial^{q} \varphi(A)}\right]
\end{array}
$$
Under the induction hypothesis, we have
$$
\begin{array}{l}
\displaystyle\frac{1}{n!}~\displaystyle \vertiii{\partial^{n} \varphi(A)}
\leq  \frac{K^{n-1}}{(2\lambda_{ min}^{{1}/{2}}(A))^3} \left[
\sum_{p+q=n-2,~p,q\geq 0}~
~\frac{C_{p}}{2^{2p}\lambda_{ min}(A)^{p}}~
~\frac{C_{q}}{2^{2q}\lambda_{ min}(A)^{q}}~
\right]
\end{array}
$$
Using the recursive formulation of the Catalan numbers $C_{k+1}=\sum_{p+q=k}~C_{p}~C_{q}$, which is valid for any $k\geq 0$ we conclude that 
\begin{eqnarray*}
\frac{1}{n!}~\vertiii{
\partial^n \varphi(A)}
&\leq& 
  \frac{K^{n-1}}{(2\lambda_{ min}^{{1}/{2}}(A))^{3}}~\frac{1}{2^{2(n-2)}\lambda_{ min}(A)^{n-2}}~
~\sum_{p+q=n-2}~C_{p}~C_{q}\\
&=&\frac{K^{n-1}}{2\lambda_{ min}^{{1}/{2}}(A)}
~\frac{C_{n-1}}{2^{2(n-1)}\lambda_{ min}(A)^{n-1}}
\end{eqnarray*}
This ends the proof of the induction. The proof of (\ref{catalan-definition}) is now completed.

We further assume that 
$$
\forall 0\leq k\leq n\qquad
\displaystyle\nabla^k \varphi(A)\cdot H
=\displaystyle \partial^k \varphi(A)\cdot H
$$
for some $n\geq 1$ and we set
$$
\Delta_n(A,H):=
\varphi(A+H)-\varphi(A)-T_n(A,H)
$$
Using (\ref{induction-square-root-integral}) we have
$$
\begin{array}{l}
\displaystyle \Delta_{n+1}(A,H)
\displaystyle=\sum_{n+2\leq k\leq 2n}\frac{1}{k!}~
\displaystyle\partial^k \varphi(A)\cdot H
-\nabla \varphi(A)\cdot \left(\left[\varphi(A+H)-\varphi(A)\right]\Delta_n(A,H)\right)\\
\\
\hskip4cm-\nabla \varphi(A)\cdot \left(\Delta_n(A,H)\left[\varphi(A+H)-\varphi(A)\right]\right)
+\nabla \varphi(A)\cdot \Delta_n(A,H)^2\\
\end{array}
$$
Under the induction hypothesis each term in the r.h.s. is of order at least $\Vert H\Vert^{n+2}$. This implies that
$$
\displaystyle\nabla^{n+1} \varphi(A)\cdot H
=\displaystyle \partial^{n+1} \varphi(A)\cdot H
$$
This yields for any $n\geq 0$ the Taylor series expansions
$$
\varphi(A+H)=\sum_{0\leq k\leq n}\frac{1}{k!}~
\displaystyle\nabla^k \varphi(A)\cdot H+\overline{\nabla}^{n+1}\varphi[A,H]$$
with the remainder term
$$
\overline{\nabla}^{n+1}\varphi[A,H]=\frac{1}{n!}~\int_0^1~(1-\epsilon)^{n}~\nabla^{n+1}\varphi\left(A+\epsilon~H\right)\cdot H ~d\epsilon
$$
To take the final step we notice that 
\begin{eqnarray*}
B=A+H&\Longrightarrow&\lambda_{ min}(A+\epsilon H)\geq(1-\epsilon)\lambda_{ min}(A)+\epsilon~\lambda_{ min}(B)\geq (1-\epsilon)\lambda_{ min}(A)
\end{eqnarray*}
This implies that
$$
\Vert\nabla^{n+1}\varphi\left(A+\epsilon~H\right)\cdot H\Vert\leq  \frac{(n+1)!}{(1-\epsilon)^{n+1/2}}
~\frac{C_{n}K^{n}}{2^{2n+1}\lambda_{ min}(A)^{n+1/2}}~\Vert H\Vert^{n+1}
$$
from which we conclude that
$$
\Vert \overline{\nabla}^{n+1}\varphi[A,H]\Vert\leq \frac{(n+1)C_{n}K^{n}}{2^{2n}\lambda_{ min}(A)^{n+1/2}}~\Vert H\Vert^{n+1}
$$
This ends the proof of the theorem.
\cqfd

\end{document}